\documentclass[11pt,oneside]{article}
\usepackage[a4paper, left=25mm, right=25mm, top=25mm, bottom=35mm]{geometry}
\usepackage{graphicx,hyperref}
\usepackage{bm,amssymb,amsfonts,amsmath,graphicx}
\newcommand{\bmb}{{\bm{b}}}
\newcommand{\bmc}{{\bm{c}}}
\newcommand{\bme}{{\bm{e}}}
\newcommand{\bmh}{{\bm{h}}}
\newcommand{\bmt}{{\bm{t}}}

\newcommand{\bmX}{{\bm{X}}}
\newcommand{\bmy}{{\bm{y}}}
\newcommand{\bmz}{{\bm{z}}}
\newcommand{\bmepsilon}{{\bm{\varepsilon}}}
\newcommand{\bmbeta}{{\bm{\beta}}}
\newcommand{\bmsigma}{{\bm{\sigma}}}

\newcommand{\matH}{{\mathbb{H}}}
\newcommand{\matI}{{\mathbb{I}}}
\newcommand{\matLambda}{{\mathbb{L}}}
\newcommand{\rootmatL}{{\mathbb{L}^{\scriptscriptstyle \! 1\! /\! 2}}}
\newcommand{\invrootmatL}{{\mathbb{L}^{\scriptscriptstyle \! -1\! /\! 2}}}
\newcommand{\matS}{{\mathbb{S}}}
\newcommand{\matX}{{\mathbb{X}}}
\newcommand{\tnorm}[1]{{\| #1 \|}}

\newcommand{\smJ}{{\!\scriptscriptstyle J}}
\newcommand{\smK}{{\!\scriptscriptstyle K}}
\newcommand{\smM}{{\!\scriptscriptstyle M}}
\newcommand{\smN}{{\!\scriptscriptstyle N}}
\newcommand{\smZ}{{\!\scriptscriptstyle Z}}
\newcommand{\smZS}{{\!\scriptscriptstyle ZS}}

\newcommand{\momy}{{\langle \bmy^2 \rangle}}
\newcommand{\momz}{{\langle \bmz^2 \rangle}}
\newcommand{\cdata}{\mathcal{D}}
\newcommand{\hmod}{\mathcal{M}}
\newcommand{\hinf}{\mathcal{H}}
\newcommand{\hmodM}{\hmod_\smM}
\newcommand{\sspace}{{\mathcal{A}}}
\newcommand{\cond}{\,|\,}

\newcommand{\bmode}{\bm{\hat{b}}}
\newcommand{\betamode}{\bm{\hat{\beta}}}

\begin{document}
\thispagestyle{empty}
\begin{center}
\textbf{\large Bayesian variable selection with spherically symmetric priors}\\[12pt]
Michiel\ B.\ De Kock and Hans\ C.\ Eggers\\[8pt]
\textit{Department of Physics, Stellenbosch University,}\\
\textit{ZA--7600 Stellenbosch, South Africa}\\
\end{center}

\begin{abstract}
  \noindent
  We propose that Bayesian variable selection for linear
  parametrisations with Gaussian iid likelihoods be based on the
  spherical symmetry of the diagonalised parameter space. Our r-prior
  results in closed forms for the evidence for four examples, including the hyper-g
  prior and the Zellner-Siow prior, which are shown to be special
  cases.  Scenarios of a single variable dispersion parameter and of
  fixed dispersion are studied, and asymptotic forms comparable to the
  traditional information criteria are derived. A simulation exercise
  shows that model comparison based on our r-prior gives good results
  comparable to or better than current model comparison schemes.
  \\

\end{abstract}

\section{Introduction}
\label{sec:introd}

The overarching problem of variable selection is to choose the best
model out of a set of candidate models $\hmodM$. Given measured data
$\cdata$, the Bayesian solution is to compute the posterior
probability for each model with Bayes' theorem,
\begin{align}
  \label{ntcb}
  p(\hmodM\cond\cdata) 
  &= \frac{p(\cdata\cond\hmodM)\,p(\hmodM)}{p(\cdata)}
  \ =\  
  \frac{p(\cdata\cond\hmodM)\,p(\hmodM)}{\sum_\smM p(\cdata\cond\hmodM)p(\hmodM)}\,.
\end{align}
As equal priors are usually assigned to the competing models, model
comparison becomes a task in finding the marginal likelihood or
evidence for each model, i.e.\ solving the integral over all $K$
model parameters $\bmbeta_\smM = (\theta_1,\ldots,\theta_\smK)$ of
the likelihood $p(\cdata\cond\bmbeta_\smM,\hmodM)$ weighted by the
parameter prior $p(\bmbeta_\smM \cond \hmodM)$,
\begin{align}
  \label{ntdb}
  p(\cdata\cond\hmodM) 
  &= \int p(\cdata\cond\bmbeta_\smM,\hmodM)\, 
  p(\bmbeta_\smM\cond\hmodM)\, d\bmbeta_\smM.
\end{align}
The preferred model will be the one with the largest evidence i.e.\
with the highest prior-weighted average over all parameters of the
likelihood.
Where computation of $p(\cdata)$ over the entire model set is
impractical or even impossible, this is circumvented by taking ratios
of two model probabilities in the form of Bayes Factors, since
$p(\cdata)$ cancels and under the equal-model-prior assumption, they
become ratios of the respective model evidences,
\begin{align}
  \label{nte}
  \mathrm{BF}(\hmodM;\hmod_{\smM'}) %
  &= \log \frac{p(\hmodM\cond\cdata)}{p(\hmod_{\smM'}\cond\cdata)}
  = \log \frac{p(\cdata\cond\hmodM)}{p(\cdata\cond\hmod_{\smM'})}\,.
\end{align}
Finding the evidence can also be difficult since model parameter
spaces $\sspace(\bmbeta_\smM)$ differ widely in size and
dimension. While convenient at first sight, assigning uniform priors
to parameters results in the untenable situation of strong dependence
of each model's evidence and consequently of Bayes Factors on
arbitrarily chosen cutoff parameters introduced by the uniform priors.
In addition, the dimension $K$ of the parameter space often differs
from model to model, compounding the problems associated with uniform
parameter priors. Furthermore, improper priors must be excluded from
the start if they are model specific because they remain in the
evidence.

These problems appear even in the simplest case of ``canonical
regression'' in which the likelihood is Gaussian and the models are
restricted to linear function spaces as studied in the past by
\cite{Jeffreys1967}, \cite{Zellner1971}, \cite{Box1973} and many
others.  %
The quest for robust and fair model comparison in this restricted
context dates back to \cite{Jeffreys1967} whose univariate Cauchy
prior was extended by \cite{Zellner1980} to multivariate form. A
simpler ``$g$-prior'' was subsequently invented by \cite{Zellner1986}
to facilitate ease of use by closed-form solutions and has found wide
application. The specific choice for $g$ and internal inconsistencies
have, however, dogged the simple $g$-prior, leading for example
\cite{Liang2008} to introduce mixtures of such $g$-priors. They showed
that $g$-mixtures resolved the inconsistencies of the simple $g$-prior
and could show it and the original Zellner-Siow prior to be special
cases within the mixture framework.

Common to all these efforts was the recognition, sometimes only
implicitly, of an underlying spherical symmetry in parameter
space. For example, the \cite{Zellner1986} prior was based on the use
of a Gaussian prior for parameters $\bmbeta$ with the same design
matrix $\matX$ as the data and precision parameter $\phi = 1/\sigma^2$
but including an additional scale parameter $g$,
\begin{align}
  \label{dce}
  p(\bmbeta\cond g,\sigma,\hinf_\smZ) 
  &= \frac{\exp\bigl[-N \bmbeta^T \matX^T\matX \bmbeta
    /2 \sigma^2 g \bigr] }
  {(\det\matX^T\matX)^{1/2}(2\pi\sigma^2 g)^{K/2}},
\end{align}
which as detailed in Section \ref{sec:nkrdiag} is easily
transformed into spherically symmetric form
\begin{align}
  \label{dcf}
  p(\bmb\cond g,\sigma,\hinf_\smZ) 
  &= \frac{e^{-N \bmb^2 / 2 \sigma^2 g}}{(2\pi\sigma^2 g)^{K/2}}.
\end{align}
As pointed out by \cite{Leamer1978}, the behaviour of the parameter
estimators is controlled by the symmetries of the prior. Often there
is no prior information which explicitly breaks the inherent spherical
symmetry of Gaussians, suggesting that spherical symmetry has been the
basis for many of the parameter priors in the canonical regression
literature all along.
 
In this paper, we take the underlying spherical symmetry to its
logical conclusion by introducing a radius variable $r$, common to all
models $\hmodM$ and for arbitrary parameter space dimension $K$, and
explicitly enforcing spherical symmetry on the hypersphere of radius $r$
by means of a $r$-prior. The projection from $\bmb$ onto
$r$ is then carried out generally, thereby reducing the $K$-dimensional
problem to a one-dimensional integral.

The $r$-prior framework introduced here encompasses earlier work as
special cases, including the conjugate-prior results of
\cite{George1997}, \cite{Raftery1997}, \cite{Berger2001}, the
various Zellner priors and the $g$-prior mixtures of \cite{Liang2008}
and shares their computational efficiency.

Not surprisingly, we find significant mathematical correspondence
between our $r$-prior and the $g$-prior mixtures of \cite{Liang2008}
as the latter implicitly assumes the same spherical symmetry made
explicit by the $r$-prior. Unlike the $g$-prior mixtures, the
$r$-prior is however not limited to mixtures of conjugate (Gaussian)
priors.

In Section \ref{sec:nkr}, we first treat the case of a single unknown
dispersion parameter $\sigma$, using it by example to introduce the
radius $r$ of the parameter hypersphere. The central result in
Eq.~(\ref{prk}) is used both to show how $g$-priors and the prior of
\cite{Zellner1980} can be obtained with particular choices of
$r$-priors as well as to introduce a simpler yet equally powerful new
$r$-prior based on properties revealed by the Mellin transform. 
In Section \ref{sec:fxr}, the single variable dispersion parameter
$\sigma$ is replaced by a set of fixed known error variances
$(\sigma_1^2,\ldots,\sigma_\smN^2)$, one for each data point. What we
have in mind here is the application of the $r$-prior formalism to
existing data with measured standard errors treating the
$\sigma_n$ not as likelihood variables but as constants.
In Section \ref{sec:bfc}, we test and compare our results to related
model comparison criteria, concluding with a discussion in Section
\ref{sec:dcn}.

\section{Single unknown dispersion parameter}
\label{sec:nkr}

\subsection{Definition and diagonalisation}
\label{sec:nkrdiag}

The generic model consists of a data set or response vector $\cdata =
\bmy = (y_1, \ldots, y_\smN) \in \mathbb{R}^N$ measured at fixed sampling
points $\bmc = (c_1, \ldots, c_\smN ) \in \mathbb{R}^N$.  The set of
predictors is represented by $K$ column vectors $\bmX_k =
(X_k(c_1),\ldots,X_k(c_\smN))^T$ which together form the $N{\times}K$
design matrix $\matX = (\bmX_1\, \bmX_2\, \cdots \, \bmX_\smK)$.
While the information $\hinf_0 = \{\bmc,N\}$ is the same for all
models, the design matrix $\matX$ and the dimensionality of the
predictor space $K$ are model-specific, $\hinf_\smM = \{\matX_\smM,
K_\smM\}$. A given model $\hmodM$ is specified by a prior $\hinf_p$ plus
$\{\hinf_0,\hinf_\smM\}$ and of course the assumption that the errors
between data and model are iid and Gaussian distributed.

From this point, we focus on developing a single model $\hmod$ and
hence drop the subscript $M$. We limit ourselves to linear regression
with coefficients $\bmbeta = (\beta_1, \ldots, \beta_\smK) \in
\sspace(\bmbeta) = \mathbb{R}^K$ and errors $\bmepsilon = \bmy -
\matX^T \bmbeta$ which are assumed to be iid and normally distributed
with a single unknown dispersion parameter, $\bmepsilon \sim
N(0,\sigma^2\bm{I_\smN})$, or
\begin{align}
  \label{ntd}
  p(\bmepsilon\cond \sigma, \hinf_0) %
  &= \prod_{n=1}^N
  \frac{e^{-\varepsilon_n^2/2\sigma^2}}{\sigma\sqrt{2\pi}},
\end{align}
resulting in the joint likelihood
\begin{align}
  \label{dgb}
  p(\bmy\cond\bmbeta,\sigma,\hmod) &= (2\pi)^{-N/2} \sigma^{-N} e^{-NQ/2\sigma^2},
\end{align}
with 
\begin{align}
  \label{dgc}
  Q(\bmy,\bmbeta,\sigma \cond \hmod) 
  &= \frac{1}{N}\tnorm{\bmy - \matX\bmbeta}^2
  = \frac{1}{N} (\bmy - \matX\bmbeta)^T (\bmy - \matX\bmbeta)
  \nonumber\\
  &= \frac{1}{N}\sum_{n=1}^N \left(y_n-{\textstyle\sum_{k=1}^K} 
    X_k(c_n)\,\beta_k\right)^2
\end{align}
related to the usual chisquared statistic by $NQ/\sigma^2 = \chi^2$.
Finding the maximum likelihood and the concomitant diagonalisation of
the parameters in $\sspace(\bmbeta)$ proceeds in the usual way, except
that we have extracted the explicit $N$-dependence in Eq.~(\ref{dgb})
and define $ \momy = \bmy^T\bmy / N$, $\matH = \matX^T\matX / N$ and
\begin{align}
  \label{dgf}
  \bmh &= \matX^T\bmy / N,
\end{align}
in terms of which
\begin{align}
  \label{dgh}
  Q &=  \momy + \bmbeta^T \matH \bmbeta - 2 \bmh^T\bmbeta.
\end{align}
The minimum of $Q$ occurs at the likelihood mode
\begin{align}
  \label{dgj}
  \betamode &= \matH^{-1}\,\bmh = (\matX^T\matX)^{-1}\matX^T \bmy.
\end{align}
The quadratic form in (\ref{dgh}) is standardised to the new parameter
set $\bmb \in \sspace(\bmb)$ via the eigenvalue equation $\matH
\,\bme_\ell = \bme_\ell\lambda_\ell$ with eigenvalues $\lambda_\ell$
and column eigenvectors $\bme_\ell$ which are orthonormalised,
$\bme^T\bme = \matI$, or using the diagonal eigenvalue matrix
$\matLambda = \mathrm{diag}(\lambda_1, \ldots, \lambda_\smN)$ and
orthogonal eigenvector matrix $\matS = (\bme_1\, \cdots\, \bme_\smK)$,
\begin{align}
  \label{dgk}
  \matH \matS &= \matS \matLambda.
\end{align}
As in \cite{Bretthorst1988}, we transform from $\bmbeta$ to $\bmb$
by a rotation by $\matS$ and a scale change by $\rootmatL =
\mathrm{diag}(\sqrt{\lambda_1}, \ldots, \sqrt{\lambda_\smK})$,
\begin{align}
  \label{dgm}
  \bmbeta &=  \matS \invrootmatL \,\bmb,\\
  \label{dgl}
  \bmb &= \rootmatL\, \matS^T \bmbeta,
\end{align}
so that the second and third terms of Eq.~(\ref{dgh})
become\footnote{Note that \cite{Bretthorst1988} uses row eigenvectors
  rather than the column vectors used in the current literature.}
\begin{align}
  \label{dgn}
  \bmbeta^T\matH\bmbeta &= \bmb^T\bmb = \bmb^2, \\
  \label{dgo}
  \bmbeta^T \bmh &= \bmb^T\bmode,\\
  \label{dgp}
  \bmode &= \rootmatL\, \matS^T \betamode,
\end{align}
and $Q$ is decomposed into a $\bmb$-independent minimum
(equivalent to minimum-$\chi^2$) and a quadratic around the mode,
\begin{align}
  \label{dgq}
  Q &= Q_0 + R_{\bmode}^2, \\
  \label{dgqb}
  Q_0 &= \tfrac{1}{N} \tnorm{\bmy - \matX\betamode}^2 
  = \momy - \betamode^T \matH \betamode
  = \momy - \bmode^2,
  \\
  \label{dgr}
  R_{\bmode}^2 &= 
  (\bmbeta - \betamode)^T \matH (\bmbeta - \betamode)
  = \tnorm{\bmb - \bmode}^2.
\end{align}
In terms of the standardised parameters, the likelihood is
\begin{align}
  p(\bmy\cond \bmb,\sigma,\hmod) %
  &=  \frac{1}{(2\pi\sigma^2)^{N/2}} 
  \exp\left[ -\frac{N}{2\sigma^2}\left(Q_0 + R_{\bmode}^2\right) \right] 
  \nonumber\\
  \label{dgz}
  &= F(\sigma) \exp
  \left[ - \frac{N}{2\sigma^2} \bmb^2 + \frac{N}{\sigma^2} \bmode^T \bmb \right],
\end{align}
with $F(\sigma) = (2\pi)^{-N/2} \sigma^{-N} e^{-N \momy /2\sigma^2}$.

\subsection{Projection onto one dimension: the $r$-prior}
\label{sec:rpr}

For model comparison, we wish to calculate the evidence, which in the present
model family is the marginal likelihood
$ %
p(\bmy\cond\hmod) 
   = \int d\bmb\,d\sigma \,p(\bmy\cond \bmb,\sigma,\hmod)
   \,p(\bmb\cond \sigma,\hmod)\,p(\sigma\cond \hmod)
$. %
Specification of the $K$-dimensio\-nal $\bmb$-prior and the integral
over $\bmb$ represents a significant challenge. In our view, the best
solution is to choose a prior for $\bmb$ which is explicitly
spherically symmetric in $\sspace(\bmb)$ by introducing a radius $r$,
\begin{align}
  \label{prc}
  p(\bmb\cond r,\hmod) %
  &= \frac{\Gamma(K/2)}{\pi^{K/2}} %
  \frac{\delta(r - \tnorm{\bmb})}{2 \,r^{K-1}} %
  \ =\ \frac{\Gamma(K/2)}{\pi^{K/2}} \frac{\delta(r^2 - \bmb^2)}{r^{K-2}},
\end{align}
where $\delta(x)$ is the Dirac delta function,\footnote{As set out in
  the literature and motivated e.g.\ by \cite{Jaynes2003appb}, the
  Dirac delta function is a limit of a sequence of probability density
  functions, and transformation of its arguments follows the standard
  rules for pdf transformation under change of variable.} plus an
intermediate $r$-prior $p(r\cond \sigma,\hmod)$. This choice of prior
is equivalent to the assumption that the prior information available
to the observer is unchanged under rotation of $\bmb$ in
$\sspace(\bmb)$.  This rotational Principle of Indifference or
``information isotropy'' in parameter space implies that $p(\bmb\cond
r,\hmod)$ must be uniformly distributed over the surface of the
$K$-dimensional hypersphere of radius $r$, for every possible value of
$r$.
Specifically, the observer has no reason a priori to prefer, or give
nonuniform prior weight to, any one of the axial directions in
$\sspace(\bmb)$ i.e.\ to any specific component of the transformed
design matrix $\matX\matS\matLambda^{-1/2}$, and hence to any original
predictor $X_k(c)$, apart from the scales and covariances introduced
by the design matrix itself during the backtransformation from $\bmb$
to $\bmbeta$.  The mathematical consequence of this argument is
Eq.~(\ref{prc}).

Once $r$ is included, the evidence is given by the $(K{+}2)$-fold
integral
\begin{align}
  \label{ntb} 
  p(\bmy\cond\hmod) %
  &= \int_0^\infty \!d\sigma\,  p(\sigma\cond\hmod) %
  \int_0^\infty \!dr \, p(r\cond\sigma,\hmod)%
  \int_{\mathbb{R}^K} \! d\bmb \, p(\bmy,\bmb \cond r,\sigma,\hmod).
\end{align}
While at first sight the extra integral may seem unnecessary, the
symmetry of prior $p(\bmb\cond r,\hmod)$ significantly simplifies the
problem since 
\begin{align}
  \label{prig}
  p(\bmy\cond r,\sigma,\hmod) %
  & = \int d\bmb\, p(\bmy,\bmb \cond r,\sigma,\hmod) 
\end{align}
can be calculated once and for all in terms of the likelihood
$p(\bmy\cond\bmb,\sigma,\hmod)$ and $r$-prior $p(\bmb\cond r,\hmod)$,
leaving us with the comparatively simple task of a two-dimensional
integral over $dr$ and $d\sigma$.

We use the Laplace-type integral representation for the Dirac delta
function and an integral representation of the generalised confluent
hypergeometric function \cite{Watson1922} 
\begin{align}
  \label{prih}
  \delta(r^2 - \bmb^2) &= \int_{\mathcal{C}}\frac{ds}{2\pi i}
  \exp\left[ sr^2 - s \bmb^2 \right], \\
  \label{prij}
  {}_0F_1(b\,;\,z) &= \frac{\Gamma(b)}{2\pi i} %
  \int_{\mathcal{C}} du\,u^{-b}\exp\left(u + \frac{z}{u}\right),
\end{align}
with $\mathcal{C}$ the contour integral along the imaginary line from
$(c{-}i\infty)$ to $(c{+}i \infty)$, to obtain 
\begin{align}
  p(\bmy\cond r,\sigma,\hmod) 
  &= \frac{F(\sigma)\,\Gamma(K/2)}{\pi^{K/2}\,r^{K-2}} 
  \int_{\mathcal{C}} \frac{ds}{2\pi i}\, e^{sr^2} 
  \int d\bmb\,\exp\left\{ - \left(\frac{N}{2\sigma^2} {+} s\right)\bmb^2 
    + \frac{N}{\sigma^2} \bmode^T \bmb
  \right\}
  \nonumber\\
  \label{prj}
  &= \frac{F(\sigma)\,\Gamma(K/2)}{\pi^{K/2}\,r^{K-2}} 
  \int_{\mathcal{C}} \frac{ds}{2\pi i} 
  \left(\frac{2\pi\sigma^2}{N{+}2\sigma^2 s}\right)^{\tfrac{K}{2}}%
  \exp\left\{ sr^2 + \frac{N^2\, \bmode^2}{2\sigma^2 (N{+}2\sigma^2 s)} \right\},
\end{align}
with $\bmode^2 = \bmh^T \matH^{-1}\bmh$ a function of $\bmy$ through
Eq.~(\ref{dgf}), leading to a closed form in terms of the generalised
hypergeometric function,
\begin{equation}
  \label{prk}
  p(\bmy\cond r,\sigma,\hmod) 
  = \frac{e^{-N(\momy+r^2)/2\sigma^2}}{(2\pi\sigma^2)^{N/2}}
  \;{}_0F_1\left(\frac{K}{2}\,;\, \frac{N^2 \bmode^2 r^2}{4\sigma^4} \right).
\end{equation}
This result is central. It shows that the sufficient statistics are
$Q_0$ and $\bmode^2$ or alternatively $\momy$ and $\bmode^2$, and
that the $K$-dimensional parameter spaces $\sspace(\bmbeta)$ and
$\sspace(\bmb)$ can be reduced to the one-dimensional space
$\sspace(r) = \mathbb{R}^+$.

The same result can be obtained via the Fourier transform
\begin{align}
  \label{prr}
  \Phi[\bmt,\bmb,p(\bmy,\bmb\cond r,\sigma,\hmod)]
  &= \int d\bmb\, e^{i\bmt^T \bmb}\,p(\bmy,\bmb \cond r,\sigma,\hmod) 
\end{align}
whose calculation proceeds exactly as above with the substitution of
$(N\bmode/\sigma^2)$ by $(N\bmode/\sigma^2) + i \bmt$, leading to
\begin{align}
  \label{prs}
  \Phi[\bmt,\bmb,p(\bmy,\bmb\cond r,\sigma,\hmod)]
  &= \frac{e^{-N(\momy+r^2)/2\sigma^2}}{(2\pi\sigma^2)^{N/2}}
  \;{}_0F_1\left(\frac{K}{2}\,;\, 
    \frac{(N \bmode + i\sigma^2 \bmt)^2 r^2}{4\sigma^4} \right),
\end{align}
from which the evidence follows as $p(\bmy\cond r,\sigma,\hmod) =
\Phi[\bmt{=}\bm{0},\bmb,p(\bmy,\bmb\cond r,\sigma,\hmod)]$.

\subsection{Connection of $r$-priors with the hyper-$g$ and Zellner-Siow priors}
\label{sec:zsg}

Before introducing a new $r$-prior, we first show that the $g$-prior
of \cite{Zellner1986}, the hyper-$g$ prior of \cite{Liang2008} and the
original Cauchy prior of \cite{Zellner1980} can all be written in
terms of suitable $r$-priors as follows.  In the case of the simple
$g$-prior, the appropriate $r$-prior is gamma-distributed,
\begin{align} \label{dcg}
  p(r\cond g,\sigma,\hinf_\smZ)  &= \frac{2}{r\Gamma(K/2)}
  \left(\frac{Nr^2}{2\sigma^2 g}\right)^{K/2}  
  e^{-N r^2 / 2\sigma^2 g},
\end{align}
leading to evidence
\begin{align}
  p(\bmy\cond g,\sigma,\hmod) 
  &= \int dr\,p(\bmy\cond r,\sigma,\hmod) \,p(r\cond g,\sigma,\hinf_\smZ) \\
  \label{ygsh}
  &= \frac{(1+g)^{-K/2}}{(2\pi\sigma^2)^{N/2}}
  \exp\left[-\frac{N\momy}{2\sigma^2}+\frac{g N\bmode^2}{2(1+g)\sigma^2}\right]
\end{align}
whose $\sigma$-integrated version can be obtained on using a Jeffreys
prior $p(\sigma\cond H_\smJ)$.

Likewise, the evidence for the hyper-$g$ prior introduced by
\cite{Liang2008}, which according to \cite{Celeux2012} is
\begin{align}
  p(\bmy\cond\hinf_g,\hmod) 
  &= \frac{(a-2)\Gamma(N/2)}{2(K+a-2)}
  \left(N \pi \momy \right)^{-N/2}
  {}_2F_1\!\left( 1\,;\, \frac{N}{2} \,;\, \frac{K{+}a}{2} \,;\, \frac{\bmode^2}{\momy}
  \right),
\end{align}
can be found either in terms of $g$ or $r$,
\begin{align}
  p(\bmy\cond\hinf_g,\hmod) 
  &= \int dr\,d\sigma\,p(\bmy\cond r,\sigma,\hmod)\,p(r\cond \sigma,K,\hinf_g)\,
  p(\sigma\cond\hinf_\smJ) \\
  &= \int dg\,d\sigma\,p(\bmy\cond g,\sigma,\hmod)\,p(g\cond \hinf_g)\,
  p(\sigma\cond\hinf_\smJ)
\end{align}
by on the one hand again using Eq.~(\ref{prk}) and a $r$-prior based
on a confluent hypergeometric function,
\begin{align}
  \label{ssd}
  p(r\cond\sigma,K,\hinf_g) 
  &= \frac{\Gamma((a{+}K)/2-1)}{\Gamma(K/2)}\, \frac{(a{-}2)}{r} 
  \left(\frac{N r^2}{2\sigma^2}\right)^{K/2}
  U\!\left(\frac{a{+}K{-}2}{2}\,;\,\frac{K}{2} \,;\, \frac{Nr^2}{2\sigma^2} \right),
\end{align}
while for the $g$-integral using Eq.~(\ref{ygsh}) and
\begin{align} \label{dci}
    p(g\cond \hinf_g) &= \frac{a-2}{2(1+g)^{a/2}},
    \quad a > 2.
\end{align}
Thirdly, the evidence for the \cite{Zellner1980} prior, which is a
complicated series of confluent hypergeometric functions
\begin{align}
  \label{ssf}
  p(\bmy\cond\hinf_{\smZS},\hmod) 
  &= 
  \sum_{j=0}^\infty \left(\frac{N \bmode^2}{2\momy}\right)^j
  \frac{\Gamma\left(\frac{1+K}{2}\right)\Gamma\left(j{+}\frac{N}{2}\right)}
  {\left(N\pi\momy\right)^{N/2} j!\,2\sqrt{\pi}}
  \;
  U\!\left(j{+}\frac{K}{2}\,;\,j{+}\frac{1}{2} \,;\, \frac{N}{2} \right),
\end{align}
can be found on the one hand in terms of $r$ using once again
Eq.~(\ref{prk}), a Jeffreys prior and a Zellner-Siow $r$-prior,
\begin{align}
  \label{ssg}
  p(r\cond\sigma,K,\hinf_{\smZS}) 
  &= \frac{\Gamma((K{+}1)/2)}{\Gamma(K/2)\,\Gamma(1/2)} %
  \;\frac{2\sigma\,r^{K-1}}{\left(\sigma^2+r^2\right)^{(1+K)/2}}.
\end{align}
Taking the alternative $g$-route by integrating
\begin{align}
  \label{ssgg}
  p(g\cond \hinf_{zs}) &= \sqrt{\frac{N}{2\pi}} e^{-N/2g} g^{-3/2}
\end{align}
together with (\ref{dcg}) and a Jeffreys prior again yields (\ref{ssf}).

\subsection{A parabolic $r$-prior} \label{sec:gmp}

Beyond the special cases covered above, the choice of $p(r\cond\hmod)$
leaves much room for new priors. In this section, we construct one
example $r$-prior, making use of the Mellin transform
\begin{align}
  \mathcal{M}(f; s) &= \int_0^\infty f(r)\,r^{s-1} dr,
\end{align}
because of its useful property of immediately exhibiting both the
asymptotic and series behaviour of the function $f(r)$. Technically,
translating the contour of the inverse Mellin transform across the
poles left of the strip of analyticity results in a series expansion
in $r$, while translation across the poles to the right gives an
asymptotic expansion. These properties are useful for examining
functions and to construct a prior with the desirable properties.

We are looking for a prior with behaviour similar to the Zellner-Siow
$r$-prior of Eq.~(\ref{ssg}) but preferably with a closed-form
solution.  The Zellner-Siow prior goes like $r^{K-1}$ close to zero
and like $r^{-2}$ for large $r$.  The Mellin transform of the
Zellner-Siow $r$-prior
\begin{align}
  \mathcal{M}(p(r\cond \sigma,K,\hinf_{\smZS}); s)
  &= \frac{\sigma^{s{-}1}}{\sqrt{\pi}}\,
  \frac{\Gamma\left[1{-}(s/2)\right]\Gamma\left[(K{+}s{-}1)/2\right]}
  {\Gamma\left[K/2\right]}
\end{align}
has a strip of convergence of $0<s<2$. Clearly,
$\Gamma\left[1-(s/2)\right]$ has poles at $s = 2,4,6,\ldots$, while
$\Gamma\left[(K{+}s{-}1)/2\right]$ has poles at $s = 1{-}K,
-1{-}K,\ldots$, which immediately gives the above desired series
expansions. This form leads, however, to a complicated evidence and so
cannot be used directly.
Taking the Mellin transform of the hyper-$g$ $r$-prior (\ref{ssd})
results in
\begin{align}
  \mathcal{M}(p(r\cond \sigma,K,\hinf_{g}); s)
  &= \frac{(a{-}2)}{2}\left(\frac{\sigma\sqrt{2}}{\sqrt{N}}\right)^{s{-}1}
  \frac{\Gamma\left[(a{-}s{-}1)/2\right]\Gamma\left[(K{+}s{-}1)/2\right]\,\Gamma\left[1{+}(s/2)\right]}
  {\Gamma\left[a/2\right]\Gamma\left[K/2\right]}.
\end{align}
The case $a=3$ is remarkably similar to the above Zellner-Siow case
and in a sense shows that the hyper-$g$ is trying to emulate the
Zellner-Siow behaviour.  Based on the above considerations, we propose
to use an $r$-prior with a similar pole structure in its Mellin
transform
\begin{align}
  \mathcal{M}(p(r\cond \sigma,K,\hinf_r); s) %
  &= \left(\frac{\sigma}{\sqrt{2 N}}\right)^{s{-}1}
  \frac{\Gamma\left[1{-}(s/2)\right]\,\Gamma\left[K{+}s{-}1\right]}{\sqrt{\pi}\,\Gamma\left[K\right]},
\end{align}
which on inversion gives us a prior in the form of a simple confluent
hypergeometric function,
\begin{align}
  p(r \cond \sigma,K,\hinf_r) 	\label{parabolic}
  &= \frac{K}{r\sqrt{\pi}}\left(\frac{Nr^2}{2\sigma^2}\right)^{K/2} %
  U\left[\frac{K{+}1}{2};\frac{1}{2};\frac{Nr^2}{2\sigma^2}\right],
\end{align}
which can also be written as a parabolic cylinder function and which
we therefore call the parabolic $r$-prior. It is of the same family as
the hyper-$g$ prior and can be reproduced by using the $g$-prior
\begin{align}
  p(g\cond\sigma,K,\hinf_r) 
  &=  \frac{\Gamma\left[1+(K/2)\right]}{\sqrt{\pi}\,\Gamma\left[\frac{1+K}{2}\right]}\frac{g^{(K-1)/2}}{(1+g)^{K/2+1}}.
\end{align}
In both cases, the resulting evidence is
\begin{align}
  p(\bmy\cond\hinf_r,\hmod)
  = \frac{\Gamma(N/2)}{2^{K+1}}\left(N \pi \momy \right)^{-N/2}
  {}_2F_1\!\left( \frac{K+1}{2}\,;\, \frac{N}{2} \,;\, K{+}1 \,;\, \frac{\bmode^2}{\momy}\right).
\end{align}
The posterior and its characteristic function are easily derived,
given the closed forms for the evidence.

\section{Known error variance}
\label{sec:fxr}

\subsection{Definition and diagonalisation}
\label{sec:fxsdiag}

In this section, we change the information from a single variable
$\sigma$ to a set of widths $\bmsigma = \{\sigma_n\}_{n=1}^N$ assumed
to be known constants, $\hinf_1 = \{\bmc,\bmsigma,N\}$, so that
the Gaussian error distribution becomes
\begin{align}
  \label{xrd}
  p(\bmepsilon\cond \hinf_1) %
  &= \prod_{n=1}^N
  \frac{e^{-\varepsilon_n^2/2\sigma_n^2}}{\sigma_n\sqrt{2\pi}}.
\end{align}
The data and predictors are now scaled individually by $\sigma_n$,
\begin{align}
  \label{xrg}
  \bmz &= \left(\frac{y_1}{\sigma_1},\ldots,\frac{y_\smN}{\sigma_\smN}\right)^T\\
  \bmX_k &= \left(\frac{X_k(c_1)}{\sigma_1},\ldots,\frac{X_k(c_\smN)}{\sigma_\smN}\right)^T\end{align}
with $\matX = (\bmX_1\,\cdots\,\bmX_\smK)$. The joint likelihood is
\begin{align}
  \label{xre}
  p(\bmy\cond\bmbeta,\hmod) 
  &= C_\sigma\,e^{-NQ/2},
\end{align}
with $C_\sigma = \bigl[\prod_n 2\pi \sigma_n^2\bigr]^{-1/2}$ a
model-independent constant and $NQ = \chi^2$ given by
\begin{align}
  \label{xrf}
  Q(\bmbeta,\bmz\cond \hmod) 
  &= \frac{1}{N}\tnorm{\bmz - \matX\bmbeta}^2
  = \frac{1}{N}\sum_n 
    \left(\frac{y_n-{\textstyle}\sum_k X_k(c_n)\beta_k}{\sigma_n}\right)^2.
\end{align}
Defining $\momz = \bmz^T \bmz / N$, $\matH = \matX^T \matX / N$ and
$\bmh = \matX^T\bmz/N$, we obtain
\begin{align}
  \label{xrk}
  Q &= \momz + \bmbeta^T \matH \bmbeta - 2 \bmh^T \bmbeta.
\end{align}
The likelihood mode in $\sspace(\bmbeta)$ is
\begin{align}
  \label{xrl}
  \betamode &= \matH^{-1}\,\bmh = (\matX^T\matX)^{-1}\matX^T \bmz, 
\end{align}
and $\bmode = \rootmatL\, \matS^T \betamode$ as before but of course
with changed $\matLambda$. Diagonalisation proceeds with the same
equations as in Section \ref{sec:nkrdiag} but subject to the above
changed definitions. We again end up with $Q = Q_0 + R_{\bmode}^2$,
with minimum
\begin{align}
  \label{xrn}
  Q_0 &= \momz - \betamode^T \matH \betamode
  = \momz - \bmode^2,
\end{align}
while $R_{\bmode}^2 = (\bmbeta - \betamode)^T \matH (\bmbeta -
\betamode)$ as before, and the likelihood itself is
\begin{align}
  p(\bmy\cond \bmb,\hmod) %
  &=  C_\sigma \exp\left[ -\frac{N}{2}\left(Q_0 + R_{\bmode}^2\right) \right] 
  \nonumber\\
  \label{xrp}
  &= C_\sigma\,e^{-N\momz/2} \exp\left[ - \frac{N}{2} \bmb^2 
    + N \bmode^T \bmb\right]
\end{align}
and the evidence for fixed $r$ changes from Eq.~(\ref{prk}) to
\begin{align}
  \label{xrc}
  p(\bmy\cond r,\hmod) &=
  C_\sigma\, e^{-N(\momz+r^2)/2}
  \;{}_0F_1\!\left(\frac{K}{2}\,;\,\frac{N^2\bmode^2 r^2}{4} \right).
\end{align}

\subsection{Results for different $r$-priors}
\label{sec:gmpf}

Since $\bmsigma$ is fixed, the parabolic $r$-prior becomes
\begin{align}
  p(r\cond K,\hinf'_r) 	
  &= \frac{K}{r\sqrt{\pi}}\left(\frac{Nr^2}{2}\right)^{K/2} U\left[\frac{K+1}{2};\frac{1}{2};\frac{Nr^2}{2}\right],
\end{align}
and the resulting evidence is
\begin{align}
  \label{xrt}
  p(\bmy\cond\hinf'_r,\hmod) 
  &= \int d\bmb\,dr\,p(\bmy\cond \bmb, \hmod)
  \,p(\bmb\cond r,\hmod)\,p(r\cond \hinf'_r)
  \\
  &= C_{\sigma} 2^{-K}e^{-N\momz/2} {}_1F_1\left(\frac{K+1}{2}\,;\,K+1\,;\, \frac{N \bmode^2}{2} \right).
  \nonumber
\end{align}
For comparison, the corresponding evidence expressions for the
hyper-$g$ and Zellner-Siow priors with their $\sigma$ set to 1 are,
respectively,
\begin{align}
  \label{xrx}
  p(\bmy\cond\hinf_g,\hmod) &= C_\sigma \, e^{-N \momz / 2} \, 
  \frac{(a-2)}{(K+a-2)}\;{}_1F_{1}\!\left( 1\,;\, \frac{K{+}a}{2} \,;\,
    \frac{N \bmode^2}{2} \right)
  \\
  \label{xry}
  p(\bmy\cond\hinf_{\smZS},\hmod)
  &= \frac{C_\sigma \, e^{-N \momz / 2}}{\sqrt{\pi}}
  \Gamma\!\left(\frac{K{+}1}{2}\right) \sum_{j=0}^\infty\frac{1}{j!}
  \left(\frac{N^2\bmode^2}{4}\right)^{\!j}
  \;U\!\left( \frac{K}{2}{+}j \,;\, \frac{1}{2}{+}j \,;\, \frac{N}{2}\right).
\end{align}
\subsection{Asymptotic forms}
\label{sec:asm}

As the argument $z$ of all the hypergeometric functions grows with
$N$, the asymptotic form for $z \gg 1$ according to
\cite{Bateman1953}
\begin{align}
  \label{asmc}
  {}_1F_{1}\!\left(a; c; z\right) &\simeq
  \frac{\Gamma(c)}{\Gamma(a)}\, z^{a-c}\, e^z
\end{align}
will often suffice. The evidence based on the parabolic $r$-prior
Eq.~(\ref{parabolic}) becomes
\begin{align}
  \label{xrsl}
  p(\bmy\cond\hinf'_r,\hmod) &\simeq \frac{C_\sigma}{\sqrt{\pi}} 
  \Gamma\!\left(\frac{K}{2}+1\right) \left(\frac{2}{N}\right)^{(K+1)/2} 
  \frac{e^{-NQ_0/2}}{\tnorm{\bmode}^{K+1}}.
\end{align}
We also find the asymptotic form of the evidence for the hyper-$g$
prior (\ref{xrx}) to be
\begin{align}
  \label{xrxl}
  p(\bmy\cond\hinf_g,\hmod) &= 
  \frac{C_\sigma  (a{-}2)}{(K{+}a{-}2)} \;
  \Gamma\!\left(\frac{K{+}a}{2}\right)
  \left(\frac{2}{N}\right)^{(K+a)/2-1}
  \frac{e^{-NQ_0/2}}{\tnorm{\bmode}^{K+a-2}},
\end{align}
and with the help of
\begin{align}
  U\!\left( \frac{K}{2}{+}j \,;\, \frac{1}{2}{+}j \,;\, \frac{N}{2}\right)
  &= \left(\frac{2}{N}\right)^{K/2+j}
  \;{}_2F_0\!\left(\frac{K}{2}{+}j \,;\, \frac{1}{2}{+}\frac{K}{2} \,;\,
  \frac{-2}{N}\right)
  \simeq \left(\frac{2}{N}\right)^{(K/2)+j},
\end{align}
approximate the Zellner-Siow evidence (\ref{xry}) by
\begin{align}
  \label{xryl}
  p(\bmy\cond\hinf_{\smZS},\hmod)
  &\simeq C_\sigma \;
  \Gamma\!\left(\frac{K{+}1}{2}\right) \left(\frac{2}{N}\right)^{K/2}
  e^{-NQ_0/2}.
\end{align}
Of course the asymptotic forms are not exactly normalised, so that we
can use them only for model comparison with information criteria or in
ratios such as Bayes Factors. 

\section{Comparing model comparison schemes}
\label{sec:bfc}

Given the closed-form expressions for the evidence within each of the
different approaches, model comparison using Bayes Factors can, of
course, be effected simply by insertion of the relevant expression
into Eq.~(\ref{nte}). We shall not do so here, however, but rather
address by example the more general question as to which of the model
comparison schemes works best. In addition to the model schemes
$\hinf_r$, $\hinf_g$ and $\hinf_{\smZS}$ considered so far, 
we include several schemes that have been used in the literature, namely %
$\hinf_{\mathrm{AIC}}$, the Akaike Information Criterion of
\cite{Akaike1974}, %
$\hinf_{\mathrm{BIC}}$, the Bayesian Information Criterion of
\cite{Schwarz1978} and %
$\hinf_{\mathrm{AICc}}$, the Akaike Information Criterion as corrected
by \cite{Hurvich1989}.  All of these can be shown to be equivalent to
$-2\,\log p(\bmy\cond\hinf)$ in our notation. For easier comparison,
we list in Table 1 the different schemes together with the
$-2\log p(\bmy\cond\hinf)$ versions of the asymptotic forms
(\ref{xrsl})--(\ref{xryl}). In the second part of Table 1, the
corresponding asymptotic forms for the evidences of Sections
\ref{sec:zsg} and \ref{sec:gmp} are shown using the relation $
{}_2F_1(a; b; c; z) \simeq \frac{\Gamma(c)}{\Gamma(a)} (bz)^{a-c}
e^{bz}$. $K$-independent constants have been omitted since they
cancel anyway once one does model comparison within any one scheme.

\begin{table}[htb]
  \begin{center}
    \begin{tabular}{|l|l|}
      \hline
      Scheme & \hspace*{15ex} $-2\log p(\bmy\cond\hinf)$ for fixed $\bmsigma$ \\ \hline \hline
      $\hinf'_r$       & $\displaystyle
      N Q_0 + \left(K+1\right)\log \!\left(\frac{N\bmode^2}{2}\right) -2\log \Gamma\!\left(\frac{K}{2}+1\right)$\\
      $\hinf_g$         & $\displaystyle
      N Q_0 + (K{+}a{-}2)\log\left(\frac{N}{2}\bmode^2\right)
      - 2 \log\Gamma\!\left(\frac{K{+}a{-}2}{2}\right)$ \\
      $\hinf_{\smZS}$     &  $\displaystyle
      N Q_0 + K\log\left(\frac{N}{2}\right) - 2 \log\Gamma\!\left(\frac{K{+}1}{2}\right)$ \\
      $\hinf_{\mathrm{AIC}}$  &  $N Q_0 + 2K$ \\
      $\hinf_{\mathrm{AICc}}$ &  $\displaystyle N Q_0 + 2 K + \frac{2K(K+1)}{N-K-1}$ \\
      $\hinf_{\mathrm{BIC}}$  &  $N Q_0 + K\log N $ \\
      \hline
      Scheme & \hspace*{15ex} $-2\log p(\bmy\cond\hinf)$ for variable $\sigma$\\ \hline \hline
      $\hinf_r$           & $\displaystyle
      -\frac{N\bmode^2}{\momy} + (K{+}1) \log \!\left(\frac{N\bmode^2}{2\momy}\right)
      -2\log \Gamma\!\left(\frac{K}{2}+1\right)$\\[6pt]
      $\hinf_g$         & $\displaystyle
      - \frac{N\bmode^2}{\momy} + (K{+}a{-}2)\log\left(\frac{N\bmode^2}{2\momy}\right)
      - 2 \log\Gamma\!\left(\frac{K{+}a{-}2}{2}\right)$ \\[6pt]
      $\hinf_{\smZS}$     &  $\displaystyle
      N\log\!\left( 1 - \frac{\bmode^2}{\momy} \right)
      + K\log\left(\frac{N}{2}\right) - 2 \log\Gamma\!\left(\frac{K{+}1}{2}\right)$ \\[6pt]
      \hline
    \end{tabular}
    \caption{Summary of model comparison schemes for the fixed
      $\bmsigma$ case of Section \ref{sec:fxr} (upper part) and for
      the variable $\sigma$ case of Section \ref{sec:nkr} (lower
      part). Constants that do not depend on $K$ are neglected.}
  \end{center}
\end{table}

In order to test our results and to make a fair comparison between
different schemes, we generate data with fixed $\bmsigma$ according to
\begin{align}
	\bmy = \matX\bmbeta + \bmepsilon 
\end{align}
where $\bmbeta$ is drawn from a Cauchy distribution centered at 0 with
its dispersion parameter set to $1$ and $5$ respectively to mimick
weak and strong signal cases. The error $\bmepsilon$ is drawn from a
standardised Gaussian distribution with a sample size $N=100$. The
design matrix is taken as orthogonal, $\matX^T\matX = \mathbb{I}_{16}$. We
use the asymptotic form of the Zellner-Siow evidence as the full form
is too slow computationally. The model size ranges successively from
$1$ to $16$ by including the first $K$ coefficients of $\bmbeta$ to
generate data $\bm{y}$ while setting the rest of the coefficients to
zero. We then calculate the highest posterior probability model using
the different priors and mean squared error loss between the fitted
and true data
\begin{align}
  \mathrm{MSE}(K) = \tnorm{\matX\bmbeta-\matX\hat{\bmbeta}^{(K)}}^2,
\end{align}
averaged over 1000 simulations.  Figure 1 shows the average MSE as a
function of model size $K$ and the model comparison schemes listed
in Table 1, including the ``Oracle'' which is the least squares
solution for the true model. To facilitate comparison, the difference
between a given method and the Oracle is shown separately in the lower
panels, while in Tables 2 and 3 the MSE values are listed for the weak
and strong signal case respectively.

We note firstly that there are large differences in the behaviour of
the model schemes for the weak and strong signal cases. At one
extreme, the BIC is quite bad for weak signals but outperforms all
other schemes for strong signals. The corrected Akaike scheme does
well for weak signals but is in mid-field for the strong signal case.
As is already apparent from the close mathematical correspondence
between the hyper-$g$ and parabolic-$r$ schemes, they converge for
large $K$ as they must. For small $K$, however, the parabolic-$r$
scheme is far superior to the hyper-$g$ scheme.

\begin{figure}
  \includegraphics[width=0.50\linewidth]{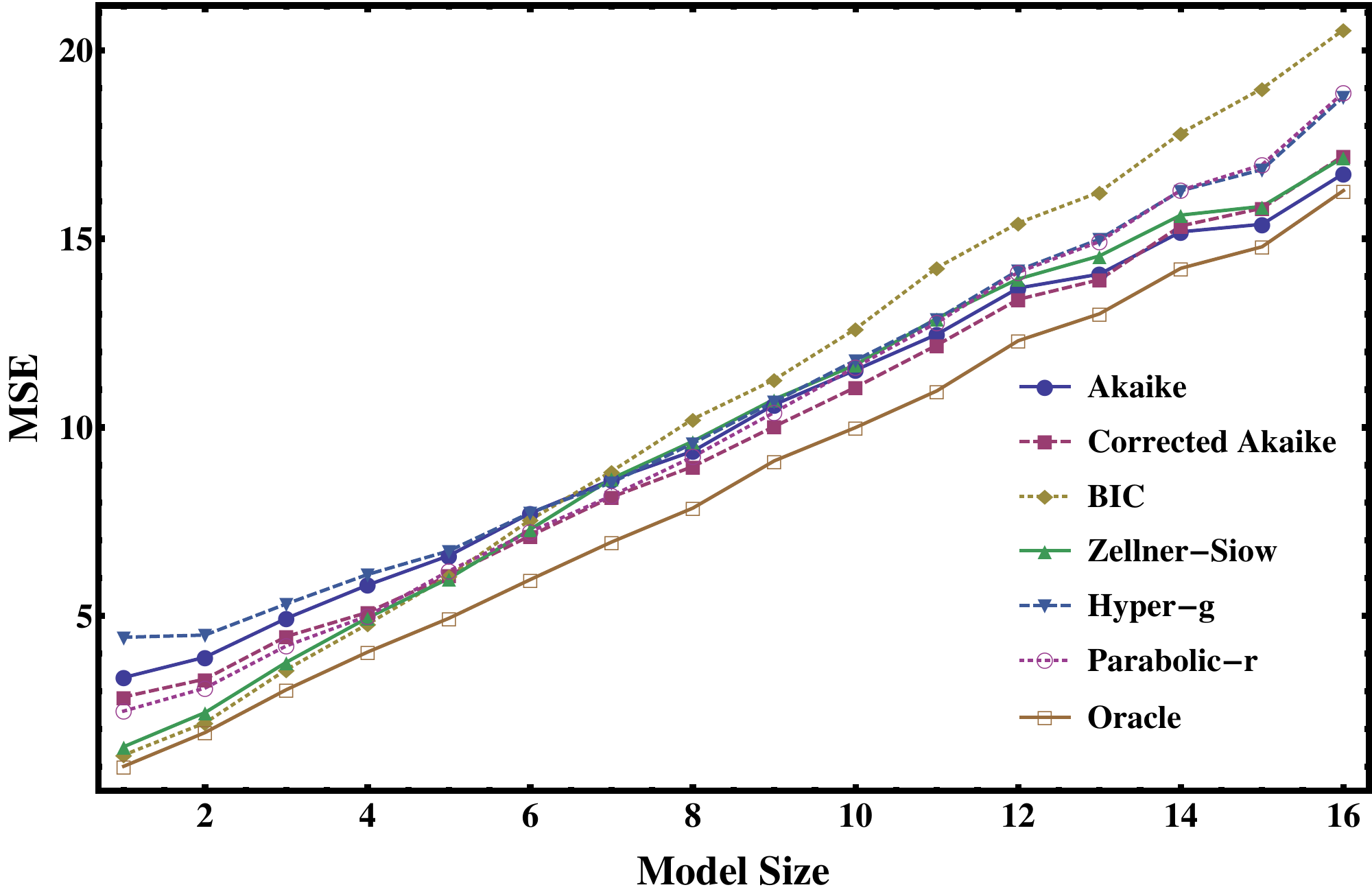}
  \includegraphics[width=0.50\linewidth]{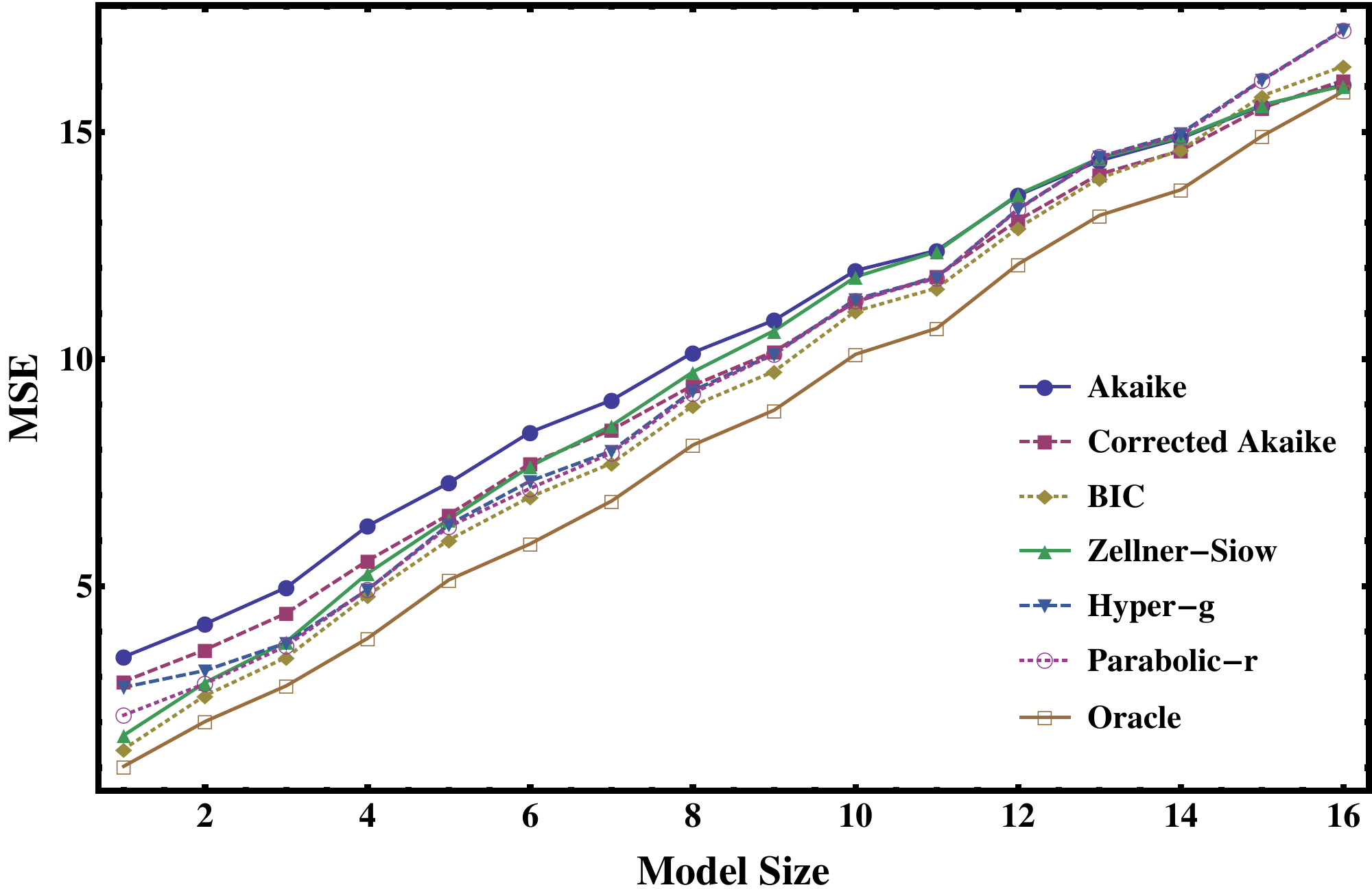}\\
  \includegraphics[width=0.50\linewidth]{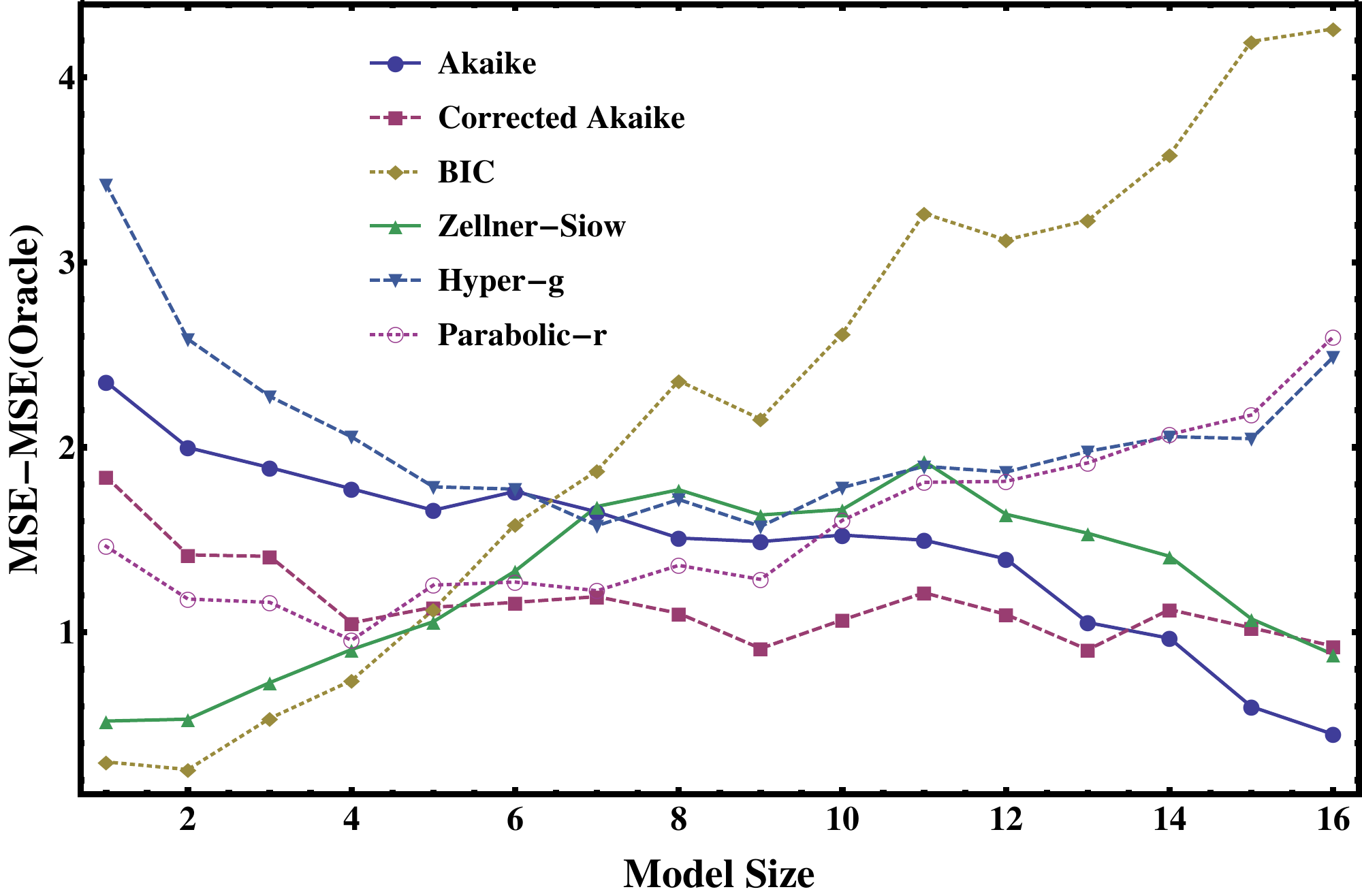}
  \includegraphics[width=0.50\linewidth]{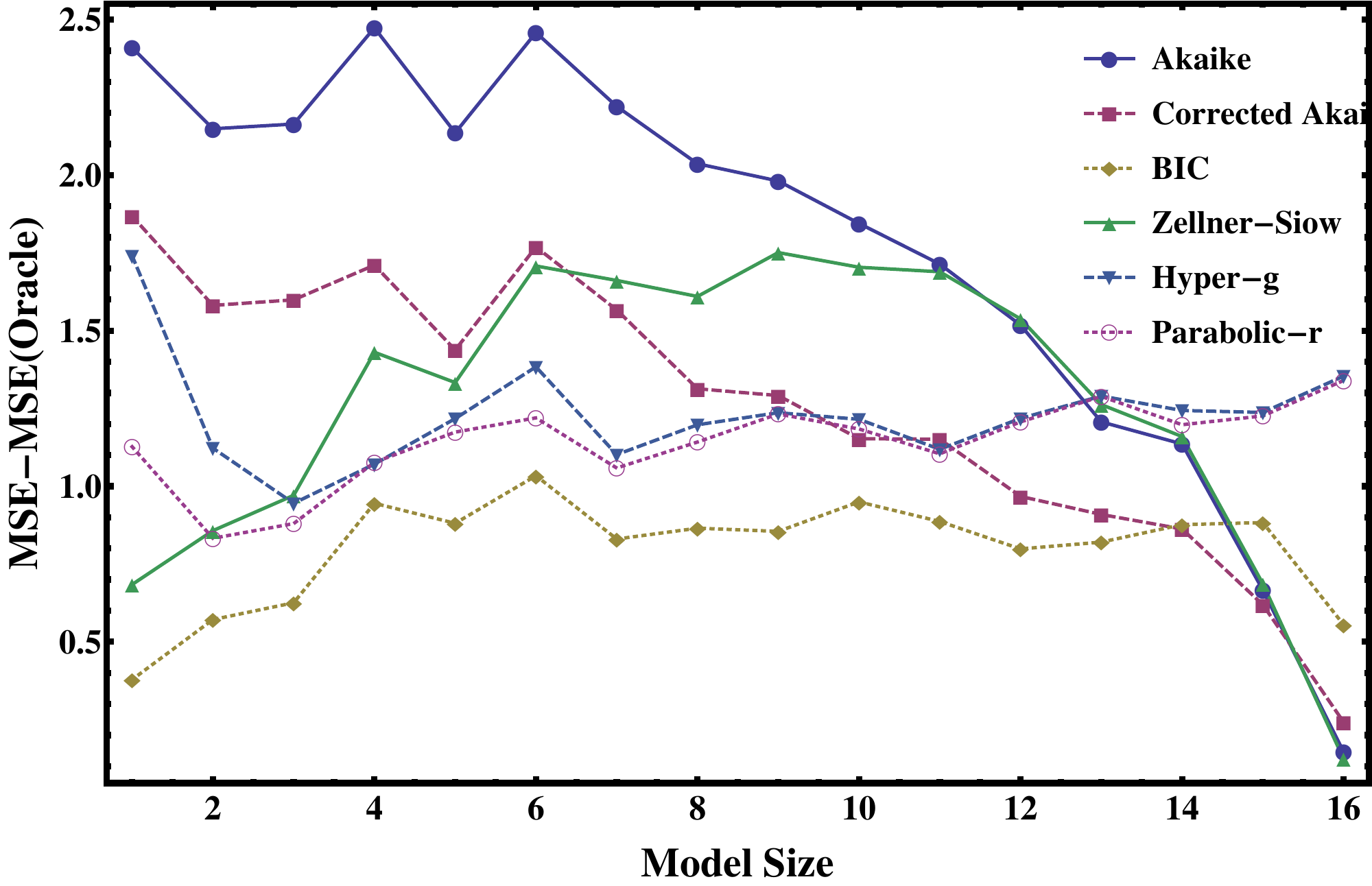}
  \caption{Upper panels: Comparison of MSE values for different model
    comparison schemes as a function of model size $K$ for weak
    signal on the left and strong signal on the right. The lower
    panels show corresponding differences between MSE($K$) and
    MSE(Oracle).}
\end{figure} 

\begin{table}[!h]
\begin{center}
  \renewcommand{\arraystretch}{1.1}
  \begin{tabular}{|c|c|c|c|c|c|c|c|}
    \hline
    $K$ & Oracle & AIC & AICc & BIC & $\hinf_\smZS$ & $\hinf_g$ & $\hinf_r$ \\ \hline \hline
    1 & 1.01 & 3.37 & 2.85 & 1.31 & 1.53 & 4.44 & 2.48 \\ 
    2 & 1.91 & 3.91 & 3.33 & 2.17 & 2.44 & 4.49 & 3.09 \\
    3 & 3.05 & 4.94 & 4.46 & 3.58 & 3.77 & 5.32 & 4.21 \\ 
    4 & 4.04 & 5.82 & 5.10 & 4.78 & 4.95 & 6.10 & 5.00 \\
    5 & 4.94 & 6.60 & 6.07 & 6.06 & 6.00 & 6.73 & 6.19 \\
    6 & 5.96 & 7.73 & 7.13 & 7.55 & 7.30 & 7.74 & 7.24 \\
    7 & 6.96 & 8.61 & 8.15 & 8.83 & 8.64 & 8.54 & 8.19 \\
    8 & 7.86 & 9.37 & 8.96 & 10.2 & 9.63 & 9.58 & 9.22 \\
    9 & 9.11 & 10.6 & 10.0 & 11.3 & 10.7 & 10.7 & 10.4 \\
   10 & 10.0 & 11.5 & 11.1 & 12.6 & 11.7 & 11.8 & 11.6 \\
   11 & 11.0 & 12.5 & 12.2 & 14.2 & 12.9 & 12.9 & 12.8 \\
   12 & 12.3 & 13.7 & 13.4 & 15.4 & 13.9 & 14.2 & 14.1 \\
   13 & 13.0 & 14.1 & 13.9 & 16.2 & 14.5 & 15.0 & 14.9 \\
   14 & 14.2 & 15.2 & 15.3 & 17.8 & 15.6 & 16.3 & 16.3 \\
   15 & 14.8 & 15.4 & 15.8 & 19.0 & 15.9 & 16.8 & 17.0 \\
   16 & 16.3 & 16.7 & 17.2 & 20.5 & 17.2 & 18.8 & 18.9 \\    \hline
  \end{tabular}
  \caption{Comparison of MSE values for different model comparison
    schemes as a function of model size $K$ for the weak signal
    case.}
\end{center}
\end{table}
\begin{table}[!h]
\begin{center}
  \renewcommand{\arraystretch}{1.1}
  \begin{tabular}{|c|c|c|c|c|c|c|c|}
    \hline
    $K$ & Oracle & AIC & AICc & BIC & $\hinf_\smZS$ & $\hinf_g$ & $\hinf_r$ \\ \hline \hline
    1 & 1.03 & 3.44 & 2.90 & 1.41 & 1.72 & 2.77 & 2.16 \\
    2 & 2.02 & 4.17 & 3.60 & 2.59 & 2.87 & 3.14 & 2.85 \\ 
    3 & 2.80 & 4.97 & 4.40 & 3.43 & 3.77 & 3.75 & 3.68 \\
    4 & 3.85 & 6.32 & 5.56 & 4.80 & 5.28 & 4.92 & 4.93 \\
    5 & 5.14 & 7.27 & 6.58 & 6.02 & 6.47 & 6.36 & 6.31 \\
    6 & 5.93 & 8.39 & 7.70 & 6.96 & 7.64 & 7.31 & 7.15 \\
    7 & 6.87 & 9.09 & 8.44 & 7.70 & 8.53 & 7.98 & 7.93 \\
    8 & 8.11 & 10.1 & 9.42 & 8.97 & 9.72 & 9.30 & 9.25 \\
    9 & 8.87 & 10.9 & 10.2 & 9.73 & 10.6 & 10.1 & 10.1 \\
   10 & 10.1 & 11.9 & 11.3 & 11.1 & 11.8 & 11.3 & 11.3 \\
   11 & 10.7 & 12.4 & 11.8 & 11.6 & 12.4 & 11.8 & 11.8 \\
   12 & 12.1 & 13.6 & 13.1 & 12.9 & 13.6 & 13.3 & 13.3 \\
   13 & 13.2 & 14.4 & 14.1 & 14.0 & 14.4 & 14.5 & 14.4 \\
   14 & 13.7 & 14.9 & 14.6 & 14.6 & 14.9 & 15.0 & 14.9 \\
   15 & 14.9 & 15.6 & 15.5 & 15.8 & 15.6 & 16.1 & 16.1 \\
   16 & 15.9 & 16.0 & 16.1 & 16.4 & 16.0 & 17.2 & 17.2 \\ \hline
  \end{tabular}
  \caption{Comparison of MSE values for different model comparison
    schemes as a function of model size $K$ for the strong signal
    case.}
\end{center}
\end{table}

\section{Discussion}
\label{sec:dcn}

We have introduced in this article the $r$-prior based on explicit
enforcement of spherical symmetry on the diagonalised parameter
space. The resulting formalism has been shown to encompass the
currently popular Zellner $g$-prior, Zellner-Siow Cauchy prior and the
hyper-$g$ prior as special cases. Beyond these, we have shown by
example of a new parabolic $r$-prior how different considerations such
as asymptotic behaviour may be incorporated. Other $r$-priors based on
further and different information can presumably be implemented in
future.

Conceptually, the $r$-priors appear to be a step towards a more formal
understanding of the symmetries on the hypersphere which are implicit
in canonical regression problems. The next step would be to understand
the scale symmetry governing $r$ itself.

The simulation shows that the $r$-prior gives good results, but also
that the detailed behaviour of it and other model schemes is quite
variable and poorly understood. Both the type of simulation and the
comparison criterion must in future be investigated in some detail.
\\

\noindent
\textbf{Acknowledgements}: This work is supported in part by a
Consolidoc fellowship of Stellenbosch University and by the National
Research Foundation of South Africa. We thank the referee for useful
comments and suggestions. Thanks also to the organisers of the
\textit{2014 ISBA--George Box Research Workshop on Frontiers of
  Statistics} for support and the participants for helpful
discussions.
\\

\bibliography{DeKockEggersP54V2}

\begin{thebibliography}{}

\bibitem[Akaike, 1974]{Akaike1974}
Akaike, H. (1974).
\newblock A new look at the statistical model identification.
\newblock {\em Automatic Control, IEEE Transactions on}, 19(6):716--723.

\bibitem[Bateman et~al., 1953]{Bateman1953}
Bateman, H., Erd{\'e}lyi, A., Magnus, W., Oberhettinger, F., and Tricomi, F.~G.
  (1953).
\newblock {\em Higher Transcendental Functions}, volume~1.
\newblock McGraw-Hill New York.

\bibitem[Berger et~al., 2001]{Berger2001}
Berger, J.~O., Pericchi, L.~R., Ghosh, J.~K., Samanta, T., and Santis, F.~D.
  (2001).
\newblock Objective bayesian methods for model selection: Introduction and
  comparison.
\newblock {\em Lecture Notes -- Monograph Series}, 38:pp. 135--207.

\bibitem[Box and Tiao, 1973]{Box1973}
Box, G. and Tiao, G. (1973).
\newblock {\em Bayesian Inference in Statistical Analysis}.
\newblock Reading, MA.

\bibitem[Bretthorst, 1988]{Bretthorst1988}
Bretthorst, G. (1988).
\newblock {\em Bayesian Spectrum Analysis and Parameter Estimation}.
\newblock Lecture notes in statistics. Springer-Verlag.

\bibitem[Celeux et~al., 2012]{Celeux2012}
Celeux, G., El~Anbari, M., Marin, J.-M., and Robert, C.~P. (2012).
\newblock Regularization in regression: Comparing bayesian and frequentist
  methods in a poorly informative situation.
\newblock {\em Bayesian Analysis}, 7(2):477--502.

\bibitem[George and McCulloch, 1997]{George1997}
George, E.~I. and McCulloch, R.~E. (1997).
\newblock Approaches for bayesian variable selection.
\newblock {\em Statistica Sinica}, 7(2):339--373.

\bibitem[Hurvich and Tsai, 1989]{Hurvich1989}
Hurvich, C.~M. and Tsai, C.-L. (1989).
\newblock Regression and time series model selection in small samples.
\newblock {\em Biometrika}, 76(2):297--307.

\bibitem[Jaynes, 2003]{Jaynes2003appb}
Jaynes, E.~T. (2003).
\newblock {\em Probability Theory: The Logic of Science (Appendix B)}.
\newblock Cambridge University Press.

\bibitem[Jeffreys, 1967]{Jeffreys1967}
Jeffreys, H. (1967).
\newblock {\em Theory of Probability}.
\newblock International Series of Monographs on Physics. Clarendon Press.

\bibitem[Leamer, 1978]{Leamer1978}
Leamer, E.~E. (1978).
\newblock Regression selection strategies and revealed priors.
\newblock {\em Journal of the American Statistical Association},
  73(363):580--587.

\bibitem[Liang et~al., 2008]{Liang2008}
Liang, F., Paulo, R., Molina, G., Clyde, M.~A., and Berger, J.~O. (2008).
\newblock Mixtures of g-priors for bayesian variable selection.
\newblock {\em Journal of the American Statistical Association}, 103(481).

\bibitem[Raftery et~al., 1997]{Raftery1997}
Raftery, A.~E., Madigan, D., and Hoeting, J.~A. (1997).
\newblock Bayesian model averaging for linear regression models.
\newblock {\em Journal of the American Statistical Association},
  92(437):179--191.

\bibitem[Schwarz, 1978]{Schwarz1978}
Schwarz, G. (1978).
\newblock Estimating the dimension of a model.
\newblock {\em The Annals of Statistics}, 6(2):461--464.

\bibitem[Watson, 1922]{Watson1922}
Watson, G. (1922).
\newblock {\em A Treatise on the Theory of Bessel Functions}.
\newblock The University Press.

\bibitem[Zellner, 1971]{Zellner1971}
Zellner, A. (1971).
\newblock {\em An Introduction to Bayesian Inference in Econometrics}.
\newblock Wiley series in probability and mathematical statistics: Applied
  probability and statistics. J. Wiley.

\bibitem[Zellner, 1986]{Zellner1986}
Zellner, A. (1986).
\newblock On assessing prior distributions and bayesian regression analysis
  with g-prior distributions.
\newblock {\em Bayesian Inference and Decision Techniques: Essays in Honor of
  Bruno De Finetti}, 6:233--243.

\bibitem[Zellner and Siow, 1980]{Zellner1980}
Zellner, A. and Siow, A. (1980).
\newblock Posterior odds ratios for selected regression hypotheses.
\newblock {\em Trabajos de estad{\'\i}stica y de investigaci{\'o}n operativa},
  31(1):585--603.

\end{thebibliography}


\end{document}